\documentclass[11pt,amssym,twoside]{article}
\usepackage{amssymb}
\newtheorem{thm}{Theorem}[section]

\newtheorem{rmk}[thm]{Remark}

\newtheorem{problem}[thm]{Problem}
\newtheorem{conj}[thm]{Conjecture}

\newtheorem{defn}[thm]{Definition}

\newtheorem{rem}[thm]{Remark}

\newcommand{\Cplx}{\mathbb C}

\newcommand{\Z}{\mathbb Z}
\newcommand{\F}{\mathbb F}
\newcommand{\G}{\mathbb G}
\newcommand{\N}{\mathbb N}
\newcommand{\Proj}{\mathbb P}
\newcommand{\Q}{\mathbb Q}
\begin{document}

\title{Self-Similarity in Geometry, Algebra and Arithmetic}%
\author{Arash Rastegar}%


\maketitle
\begin{abstract}
We define the concept of self-similarity of an object by
considering endomorphisms of the object as `similarity' maps. A
variety of interesting examples of self-similar objects in
geometry, algebra and arithmetic are introduced. Self-similar
objects provide a framework in which, one can unite some results
and conjectures in different mathematical frameworks. In some
general situations, one can define a well-behaved notion of
dimension for self-similar objects. Morphisms between
self-similar objects are also defined and a categorical treatment
of this concept is provided. We conclude by some philosophical
remarks.
\end{abstract}
\section*{Introduction}

The goal of this manuscript is to study the concept of
self-similarity in the abstract sense of the word. For example, an
abstract notion of equality can be implemented by fixing an
equivalence relation. In the same way, given any concept of
similarity, we consider the associated concept of a self-similar
object. For now, consider the following general definition for
self-similar objects:

\begin{defn}
Let $X$ be a space and let $f_i$ for $ i=1$ to $n$ denote
appropriate endomorphisms of $X$. A subset $F\subseteq X$ is
called `self-similar with respect to $f_i$' if $F$ is disjoint
union of its images under the endomorphisms $f_i$. We write
$F=\sqcup_i f_i(F)$. Endomorphisms $f_i$ are called `similarity
maps'. $F$ is called `self-similar' if it is self-similar with
respect to a finite number of similarity maps. Evidently, choice
of $f_i$ is not unique. If finite intersection of similar images
are allowed, $F$ is called an `almost self-similar' subset of $X$.
\end{defn}

This formulation works for both geometric and arithmetic settings
and is even fit for an algebraic treatment after minor changes.
Several examples are introduced in the first section, and many old
results and conjectures are reformulated using this language in
the second section. A more abstract categorical formulation paves
the way for a functorial understanding of self-similarity which
form the content of the third chapter. Philosophical remarks
conclude this manuscript.

\section{Examples of self-similar objects}

In order to provide grounds for better understanding of the scope
of generality of results proposed in the following sections, we
start by providing a list of self-similar objects. In geometry,
there is a rich gallery of spaces and endomorphisms. Some
classical self-similar affine fractals are already introduced in
the literature. In algebra, considering automorphisms of
self-similar objects and associated algebraic objects lead to
some self-similar objects. In algebraic geometry, algebraic
endomorphisms of varieties could be considered as similarity
maps. In arithmetic, there are fairly general frame-works in
which we have an understanding of the asymptotic behavior of
points of bounded norm or bounded height. We will see that, these
are natural places to look for self-similar arithmetic objects.
Some of these arithmetic self-similar objects can be formulated in
the language of arithmetic geometry.

\subsection{Self-similarity in geometry}

Thinking of self-similar geometric objects as geometric spaces
would lead to an unfamiliar concept of geometric intuition. Here
are some examples:

$\bold 1.$ Cantor set, Serpinski carpet, von Koch curve are all
examples of affine self-similar or almost self-similar objects. To
obtain such objects, one starts from a simple geometric object and
finitely many non-intersecting similar copies inside. By iteration
of this process one always gets a self-similar object with
respect to similarity maps one starts with. To obtain Serpinski
carpet and von Koch curve one allows one point of intersection.
These are examples of affine self-similar objects, for which one
starts from an affine ambient space and considers finitely many
affine maps as similarity maps. Standard box formulas provide a
notion of fractal dimension [Falc].

$\bold 2.$ Take the projective space $\Proj^n(\Cplx )$ as ambient
space and projective endomorphisms as similarity maps. This way
one could embed affine fractals in projective ambient spaces.
There are examples of projective self-similar objects which are
not induced by affine objects. For example, the subset
$$
\{(2^i;2^j)\in \Proj^1(\Cplx)|i,j\in \N \cup \{ 0\}\}\cup \{ (0;1)
,(1;0)\}
$$
is a projective self-similar set with $f(x;y)=(x;2y)$ as a single
similarity map.

$\bold 3.$ The set of periodic points and the set of pre-periodic
points of an endomorphism of a geometric space could also be
thought as a self-similar object. Kawaguchi proved that for a
holomorphic endomorphism of degree $\geq 2$ of the projective
space $\Proj^n(\Cplx)$ or of a compact Riemann surface there are
only finitely many periodic points [Ka]. This brings finite
dynamical systems into attention. There are many other finiteness
results in this direction. For example, Fakhruddin proved that
the set of all possible periodicity exponents for an endomorphism
of a proper projective variety is a finite set [Fakh].

\subsection{Self-similarity in algebra}

Algebraic objects associated to self-similar geometric objects
usually express a self-similar nature:

$\bold 1.$ Consider the vector space $V$ of continuous real-valued
functions on a geometric self-similar object with injective
similarity maps. Any similarity map induces a self-map of this
vector space if we extend by zero. The vector space $V$ is
generated by finitely many non-intersecting copies of its similar
images which makes a self-similar vector space.

$\bold 2.$ Take a ring $R$ as ambient space and injective
endomorphisms as similarity maps. Suppose that $R$ can be
generated by finitely many of its own copies. This makes $R$ to
be a self similar object. For example, $\Z$ can be generated by
finitely many subrings isomorphic to itself.

$\bold 3.$ Consider the Cayley graph of a finitely generated free
group with respect to a free generating set. Remove the identity
vertex and edges landing on it. Let $G$ denote the automorphism
group of one of the connected components of the remaining graph.
$G$ is generated by a finite number of subgroups isomorphic to
$G$ together with finitely many of their cosets. Therefore, $G$
is a self-similar object.

\subsection{Self-similarity in algebraic geometry}

In algebraic geometry, an algebraic variety plays the role of the
ambient space and algebraic endomorphisms of varieties could be
taken as similarity maps:

$\bold 1.$  Take the complex projective space $\Proj^n(\Cplx )$
as ambient space and polynomial endomorphisms
$f_i=(\phi_0,\phi_1,...,\phi_n)$ which are induced by $n+1$
homogeneous polynomials in $n+1$ variables of degree $m_i$, as
similarity maps. For example, the subset
$$
\{(2^i;2^j)\in \Proj^1(\Cplx)|i,j\in \N \cup \{ 0\}\}
$$
is self-similar with respect to $f_1(x;y)=(x^2;y^2)$ and
$f_2(x;y)=(2x^2;y^2)$.

$\bold 2.$ Let $A$ denote an abelian variety over $\Cplx$. Fix a
natural number $n$. Take $A$ for an ambient space and translated
multiplications by $n$ as similarity endomorphisms. Then, for
every integer $n>1$, the set of torsion elements of $A$ is a
self-similar subset with respect to all translations of
multiplication by $n$ via $n$-torsion elements.

$\bold 3.$ Let $A$ denote an abelian variety over $\Cplx$. A
finitely generated subgroup of $A$ is self-similar with respect
to a few translations of multiplication by $n$ taken as
similarity maps if $n$ is appropriately chosen. The following
result of Neron addresses the growth of the number of points of
bounded height [Ne]:

Let $A\subset \Proj^n$ denote an abelian variety defined over a
number field  $K$ and let $r=r(A,K)$ denote the rank of the group
of $K$-rational points in $A$, then there exists a constant
$c_{A,K}$ such that
$$
N(A(K),x)\sim c_{A,K}(log x)^{r/2}.
$$
where $N(A(K),x)$ denotes the number of points of height bounded
by $x$. One can see that the rate of growth mentioned above is
related to fractal dimension of the set of $K$-rational points as
a self-similar subset of $A$ [Ra].

\subsection{Self-similarity in arithmetic}

The setting of arithmetic self-similarity allows us to define a
well-defined and well-behaved concept of fractal dimension for
these arithmetic objects. This concept is intimately related to
the concept of arithmetic height [Ra]:

$\bold 1.$ Take a ring $R$ as ambient space and polynomial maps as
as similarity maps. A self-similar subset is affine, in case
these polynomials are of degree one. An example of affine
self-similar subset of $\Z$ is the set of integers with only
certain digits appearing in their decimal expansion. The fractal
dimension of an affine self-similar subset of $\Z$ is defined by
the box formula $\sum_i a_i^{-s}=1$ where $a_i$ are the leading
coefficients of degree one similarity maps. This is a
well-defined and well-behaved notion of fractal dimension [Ra].

$\bold 2.$ By a $t$-module of dimension $N$ and rank $d$ defined
over the algebraic closure $\bar k=\overline {\F_q(t)}$ we mean,
fixing an additive group $(\G_a)^N(\bar k)$ and an injective
homomorphism $\Phi$ from the ring $\F_q[t]$ to the endomorphism
ring of $(\G_a)^N$ which satisfies
$$
\Phi(t)=a_0F^0+...+a_dF^d
$$
where $a_i$ are $N\times N$ matrices with coefficients in $\bar
k$ with $a_d$ non-zero, and $F$ is a Frobenius endomorphism on
$(\G_a)^N$. Polynomials $P_i\in \F_q[t]$ of degrees $r_i$ for $
i=1$ to $n$ could be taken as similarity maps. Denis defines a
canonical height $\hat h$ on these modules for which one has
$$
\hat h[\Phi(P)(\alpha)]=q^{dr}.\hat h[\alpha]
$$
for all $\alpha \in (\G_a)^N$, where $P$ is a polynomial in
$\F_q[t]$ and degree $r$ [De]. We define the fractal dimension of
a self-similar subset with respect to $P_i$ to be the real number
$s$ such that $\sum_i (r_id)^{-s}=1$.

$\bold 3.$ Start from a linear semi-simple algebraic group $G$
and a rational representation $\rho :G\to GL(W_{\Q})$ defined over
$\Q$. Let $w_0\in W_{\Q}$ be a point whose orbit $V=w_0\rho (G)$
is Zariski closed. Then the stabilizer $H\subset G$ of $w_0$ is
reductive and $V$ is isomorphic to $H\setminus G$. By a theorem of
Borel-Harish-Chandri $V(\Z)$ breaks up to finitely many $G(\Z)$
orbits [Bo-HC]. Thus the points of $V(\Z)$ are parametrized by
cosets of $G(\Z)$. Fix an orbit $w_0G(\Z)$ with $w_0$ in $G(\Z)$.
Then the stabilizer of $w_0$ is $H(\Z)=H\cap G(\Z)$.
Duke-Rudnick-Sarnak [D-R-S] putting some extra technical
assumptions, have determined the asymptotic behavior of
$$
N(V(\Z),x)=\circ\{\gamma \in H(\Z)\setminus G(\Z):
||w_0\gamma||\leq x\}.
$$
They prove that there are constants $a\geq0 ,b>0$ and $c>0$ such
that
$$
N(V(\Z),x)\sim cx^a(log x)^b.
$$
The additive structure of $G$ allows one to define self-similar
subsets of $V(\Z)$ and study their asymptotic behavior using the
idea of fractal dimension. The whole set $V(\Z)$ could not be a
fractal of finite dimension, since the asymptotic behavior of its
points is not polynomial [Ra].

\subsection{Self-similarity in arithmetic geometry}

In arithmetic geometry, the geometric nature of self-similar
objects is combined with arithmetic and height function formalism:

$\bold 1.$ Take the complex projective space $\Proj^n(\Cplx )$ as
ambient space and polynomial endomorphisms
$f_i=(\phi_0,\phi_1,...,\phi_n)$ which are induced by $n+1$
homogeneous polynomials in $n+1$ variables of degree $m_i$, as
similarity maps. $\Proj^n(K)$ is self-similar for any given
number field $K$, but with respect to infinitely many similarity
maps. Schanuel describes the asymptotic behavior of points in
$\Proj^n(K)$ [Scha]: Let $h,R,w,r_1 ,r_2 ,d_K ,\zeta_K$ denote
class number, regulator, number of roots of unity, number of real
and complex embeddings, absolute discriminant and the zeta
function associated to the number field $K$. Then the asymptotic
behavior of points in $\Proj^n(K)$ of non-logarithmic height
bounded by $x$ is given by
$$
{hR \over {w\zeta_K(n+1)}}\left( {2^{r_1}(2\pi)^{r_2} \over
{d_K^{1/2}}}\right)^{n+1}(n+1)^{r_1+r_2-1}x^{n+1}.
$$
One can show that fractal dimension of $\Proj^n(K)$ is related to
the above asymptotic behavior which implies that $\Proj^n(K)$ is
not a fractal.

$\bold 2.$ Schmidt [Schm] (in case $K=\Q$) and Thunder [Th] for
general number field $K$ generalized the above theorem to
Grassmanian varieties, and proved that
$$
C(G(m,n)(K),x)\sim c_{m,n,K}x^n
$$
where $C$ denotes the number of points of bounded non-logarithmic
height and $c_{m,n,K}$ is an explicitly given constant. This
implies that $G(m,n)(K)$ is not a fractal either. The following is
a generalization to flag manifolds proved by
Franke-Manin-Tschinkel [Fr-Ma-Tsh]):

Let $G$ be a semi-simple algebraic group over $K$ and $P$ a
parabolic subgroup and $V=P\backslash G$ the associated flag
manifold. Choose an embedding of $V\subset\Proj^n$ such that the
hyperplane section $H$ is linearly equivalent to $-sK_V$ for some
positive integer $s$, then there exists an integer $t\geq 0$ and
a constant $c_V$ such that
$$
C(V(K),x)^s=c_Vx(log x)^t.
$$

$\bold 3.$ Let $\F_q(X)$ denote the function field of an
absolutely irreducible projective variety $X$ which is
non-singular in codimension one, defined over a finite field
$\F_q$ of characteristic $p$. one can consider the logarithmic
height on $\Proj^n(\F_q(X))$ defined by Neron [La-Ne] and study
the asymptotic behavior of $N(\Proj^n(\F_q(X)),d)$. Serre and Wan
proved the following [Wa]: Let $n(\Proj^n(\F_q(X)),d)$ denote the
number of points in $\Proj^n(\F_q(X))$ with logarithmic height
equal $d$ then
$$
n(\Proj^n(\F_q(X)),d)={{hq^{(n+1)(1-g)}} \over
{\zeta_X(n+1)(q-1)}}q^{(n+1)d}+O(q^{d/2+\epsilon})
$$
for $d\to \infty$. Therefore $N(\Proj^n(\F_q(X)),d)\sim
cq^{(n+1)d}$ for a constant $c$, which confirms that
$\Proj^n(\F_q(X))$ could be thought of as a fractal.

\section{Unification results and conjectures}

The general formulation of self-similarity introduces a new
perspective to many old results. Many seemingly unrelated results
can be unified in a single self-similar framework. In this
section, we provide perspectives and formulate some results and
conjectures in this new setting.

\subsection{Prospects of self-similarity in geometry}

Let us start with self-similar vector spaces. Algebraic
operations on vector spaces such as direct sum, tensor product,
wedge and symmetric powers induce new self-similar vector spaces.
The concept of linear endomorphism, should be restricted to
linear maps commuting with all similarity maps. To define a
morphism between different self-similar vector spaces, one should
associate to any similarity map of the origin a similarity map of
the target such that these two commute the linear map between
self-similar vector spaces. This gives us a perspective towards a
formulation of self-similarity in the categorical setting.

The concept of pairing and duality and other related concepts of
linear algebra could be extended to the world of self-similar
vector spaces. As usual, pairings take values in the base field
and are not supposed to be compatible with similarity maps. Dual
of a self-similar vector space is again self-similar and fits in
a nondegenerate pairing with the original vector space.

Suppose that a self-similar vector space can be generated by
finitely many vectors and their images under all combinations of
similarity maps. This introduces the notions of basis and
dimension, and brings us to the realm of matrix theory. Although
one can not determine a self-similar linear map between
self-similar vector spaces with finitely many numbers, a
beautiful theory of determinants survives the new complications.

Self-similar vector spaces introduce a new perspective to the
concept of space. The new formulation of linear algebra brings up
a whole new world of geometrical objects. Self-similar manifolds
could be defined as spaces which are locally isomorphic to
self-similar open subspaces of self-similar vector spaces.
Self-similar tangent spaces and self-similar vector bundles and
the corresponding geometric set up could be studied. Self-similar
cohomology theories could also be introduced which satisfy
Poincare duality. These cohomology groups are self-similar vector
spaces with similarity maps induced by the similarity maps of the
corresponding space. Self-similar intersection theory and
self-similar $K$-theory could also be studied in this new setting.

Self-similar spaces could be used as mathematical models for the
flow of time and continuous movement. This brings up a new
perspective towards space-time and possibility to formulate the
infinitesimal notions of quantum theory in the new geometric
setting. In fact, real line which appeared first as a model for
the flow of time is the origin of the Euclidean concept of
similarity.

\subsection{Prospects of self-similarity in algebra}

Since a module is nothing but a family of vector spaces, the
concept of a self-similar vector space, when considered in
families, defines a self-similar module. A self-similar module is
generated by finitely many non-intersecting submodules each
isomorphic to the full module. A morphism between self-similar
modules is defined similar to a morphism between self-similar
vector spaces. Direct sum, tensor product, wedge and symmetric
powers of self-similar modules induce new self-similar modules.

Main concepts of homological algebra could also be formulated in
the new setting of self-similar modules. The main obstacle is
that a self-similar module can not be free. We shall replace the
concept of self-similar free resolution by almost self-similar
free resolution where an almost self-similar free module with
respect to $n$ similarity maps $\phi_1,..., \phi_n$ is defined as
the direct sum $F(n)$ of a free module $F$ and all formal copies
$\phi_{i_1}\circ ...\circ \phi_{i_k}(F)$. It is evident how to
define self-similarity maps $\phi_i:F(n)\to F(n)$.

The automorphism groups of finite dimensional self-similar vector
spaces which are self-similar analogue of classical groups could
also be studied. The automorphism group of a self-similar vector
space is a self-similar group. It is a beautiful challenge to
develop the whole theory of Lie groups and Lie algebras and their
representations in the new self-similar setting. A self-similar
representation is nothing but an action of a self-similar object
on a self-similar vector space.

The concept of self-similar classical groups brings up the
concept of self similar Hopf algebras. This motivates how to
define a self-similar ring and a self-similar algebra. A
self-similar ring is a ring generated as a ring by finitely many
subrings each isomorphic to itself. Self-similar rings and
self-similar modules are intimately related. Self-similar ideals
induce self-similar rings after taking the quotient. Morphisms
between self-similar rings is defined as usual.

Let a group act on a self-similar object by self-similar
endomorphisms. The fixed subring will also be self-similar. One
could even consider the action of a self-similar group on a
self-similar ring. To get an idea how to define a self-similar
action, one shall consider the action of the automorphism group
$G$ of a self-similar vector space $V$ on the vector space
itself. To each self-map $\psi_i$ of $G$ one shall associate a
self-map $\phi_j$ of $V$ such that for all $g\in G$ and $v\in V$
we have $\phi_j(g.v)=\psi_i(g).\phi_j(v)$. The corresponding
invariant theory could also be developed in the setting of
self-similar actions.

Self-similar deformation of self-similar algebraic structures is
also an important topic. This makes a bridge between self-similar
geometric spaces and self-similar algebraic structures. One
should preserve the self-similarity structure while deforming,
but choice of self-similar maps is not unique. This is why the
deformation space of a self-similar algebraic object is a
self-similar space.

\subsection{Prospects of self-similarity in algebraic geometry}

We get fractal versions of Siegel's theorem on finiteness of integral points
and Falting's theorem on diophantine approximation on abelian varieties. 
Here are special cases, which could be formulated without any reference to fractals.
We have treated these special cases separatedly in [Ras1] and [Ras2]:

\begin{thm}
Let $X$ be an affine open subcurve of a connected smooth projectuve curve of genus $\geq 1$ defined
over $\mathbb{C}$  in the ambient affine space $\mathbb{A}^n(\mathbb{C})$ and let $F\subset \mathbb{A}^n(\mathbb{C})$
denote any finitely generated subgroup of $\mathbb{C}^n$ . Then $X(K )\cap F$ is finite.
\end{thm}

This implies that Siegel's theorem is an algebro-geometric fact, not an Arithmetic one.

\begin{thm}
Let $A$ be an abelian variety defined over a finitely generated subfield $K$ of $\mathbb{C}$. Let
$E$ is a geometrically irreducible subvariety of $A$ defined over $K$ and $F$ be a finitely generated subgroup of 
$A(K)$.                                                                                                                                                                                                                                                                       
Let $w$ be a valuation on $K$ and $H(x)$ a height function on $K$ coming from a choice of projective model for $K$
over the algebraic closure of $\mathbb{Q}$ in $K$.
If $d_w(x,E)$ denotes the $w$-adic distance from $x$ to $E$, and $\kappa$ and $c$ are positive constants,
then, there are only finitely many points in $F$ satisfying the following inequality
$$
d_w(x,E)< cH(x)^{-\kappa}.
$$
\end{thm}

This, in turn, implies that Faltings' theorem on Diophantine approximation on abelian varieties is also 
an algebro-geometric fact, not an Arithmetic one.

Arithmetic fractals provide a common framework in which similar
theorems in Diophantine geometry could be united in a single
context. Considering the fact that one could think of $A(\bar K)_{tor}$
and $A(K)$ as fractals in $A$, the theorems of Raynaud [Ra] and
Faltings [Fa] can be united in the following general conjecture:

\begin{conj} (Fractal conjecture)
Let $V$ be an irreducible variety defined over a finitely
generated field $K$ and let $F\subset V(K)$ denote a fractal
on $V$ with respect to finitely many height-increasing self-maps
$$
f_i:V \rightarrow V 
$$ 
defined over $K$.
Then, for any reduced subscheme $Z$ of $V$ defined over
$K$ the Zariski closure of $Z(\bar K)\cap F$ is union of finitely
many points and finitely many components $B_j$ such that $B_j(K)\cap F$ 
is a fractal in $B_j$ for each $j$, with respect to some of $f_i$.
If non of the components of $B_j$ are pre-priodic with respect to any of $f_i$
then any $B_j(K)\cap F$ 
is a fractal in $B_j$ with respect to all of $f_i$.
\end{conj}

\begin{rmk}
You can start with $F\subset V(\bar K)$, but then you can not assume 
$f_i$ are height increasing and instead you may join some $B_j$ to make a fractal.
\end{rmk}

In particular, we have stated the following 
 
\begin{conj}
For any algebraic curve $C$ embedded
in $V$ defined over $K$ which is not invarient under $f_i$, 
we have $C(\bar K)\cap F$ is at most a finite set.
\end{conj}

A consequence of the above conjecture would be the generalized
Lang's conjecture which is confirmed by results of Raynaud [Ray],
Laurent [Lau], Zhang [Zh]:

\begin{conj}
(Lang) Let $X$ be an algebraic variety defined over a number-field
$K$ and let $f:X \to X$ be a surjective endomorphism defined over
$K$. Suppose that the subvariety $Y$ of $X$ is not pre-periodic in
the sense that the orbit $\{Y,f(Y),f^2(Y),...\}$ is not finite,
then the set of pre-periodic points in $Y$ is not Zariski-dense
in $Y$.
\end{conj}

Here is another implication of our self-similarity conjecture in the spirit of above conjecture by Lang:

\begin{conj}
(Forward orbit conjecture) Let $V$ be an irreducible variety defined over a finitely
generated field $K$ and let $f_i:V \rightarrow V$ denote finitely many self maps of $V$ defined over $K$. 
Let $F$ denote the forward orbit with respect to $f_i$ of finitely many points of $V$ defined over $K$. 
Then, for any reduced subscheme $Z$ of $V$ defined over
$K$ the Zariski closure of $Z(\bar K)\cap F$ is union of finitely
many points and finitely many components $B_j$ such that $B_j(\bar K)\cap F$ 
is the forward orbit with respect to some $f_i$ of finitely many points of $B_j$ defined over $K$, for each $j$.
\end{conj}

It is instructive to notice that, the common geometric structures
appearing in the context of Diophantine geometry, is exactly the
same as the objects appearing in dynamics of endomorphisms of
algebraic varieties which was the original context that height
functions were introduced.

We also present another conjecture in the same lines for
quasi-fractals in an algebraic variety $X$, where
self-similarities are allowed to be induced by geometric
self-correspondences on $X$ instead of self-maps. This time, we
drop the requirement that similar images shall be almost-disjoint.

\begin{conj}
(Quasi-fractal conjecture)
Let $V$ be an irreducible variety defined over a finitely
generated field $K$ and let $F\subset V(\bar K)$ denote a
quasi-fractal on $V$ with respect to correspondences $Y_1,...,Y_n$
on $V$ living in $V\times V$ with both projections finite and
surjective. $F$ may contain
a subvariety of $V$. Then, for any reduced subscheme $Z$ of $V$ defined
over $K$ the Zariski closure of $Z(\bar K)\cap F$ is union of
finitely many points and finitely many components $B_i$ such that
for each $i$ the intersection $B_i(\bar K)\cap F$ is a
quasi-fractal in $B_i$ with respect to some correspondences induced by
$Y_i$.
\end{conj}

The above conjecture is confirmed by Andre-Oort
conjecture which are proved today.
\subsection{Prospects of self-similarity in arithmetic}

The idea of considering self-similar subsets of $\Z$ is due to O.
Naghshineh who proposed the following problem for "International
Mathematics Olympiad" held in Scotland in July 2002.
\begin{problem}
Let $F$ be an infinite subset of $\Z$ such that $F=\bigcup_{i=1}^n
a_i.F+b_i$ for integers $a_i$ and $b_i$ where $a_i.F+b_i$ and
$a_j.F+b_j$ are disjoint for $i\neq j$ and $|a_i|>1$ for each $i$.
Prove that
$$
\sum_{i=1}^n {1\over |a_i|}\leq 1.
$$
\end{problem}

In [Na], he explains his ideas about self-similar sunsets of $\Z$
and suggests how to define their dimension:

\begin{defn}
Let $\phi_i:\Z \to \Z$ for $ i=1$ to $n$ denote linear maps of
the form $\phi_i(x)=a_i.x+b_i$ where $a_i$ and $b_i$ are integers
with $|a_i|>1$. Let $F\subseteq \Z$ be a fractal with respect to
$\phi_i$. The fractal dimension of $F$ is defined to be the real
number $s$ such that
$$
\sum_{i=1}^n |a_i|^{-s}=1.
$$
\end{defn}

\begin{thm} (Mahdavifar, Naghshineh) Let $F_1\subseteq F_2 \subseteq \Z$
be fractals. Then the notion of fractal dimension is well-defined
and $dim(F_1) \leq dim(F_2)$.
\end{thm}

The proof of the above theorem works for $\Z[i]$ as well. This
motivates us to consider the general case of arbitrary number
field. Let $K$ be a number field and let $O_K$ denote its ring of
integers. One can take $O_K$ as ambient space and polynomial maps
$\phi_i:O_K \to O_K$ of degrees $n_i$ with coefficients in $O_K$
as self-similarities. Let $a_i$ denote the leading coefficient of
$\phi_i$, and $n_i$ denote the degree of $\phi_i$. Fix an
embedding $\rho :K\hookrightarrow \mathbb C$. Assume
$Norm(a_i)>1$ in case $\phi_i$ is linear. Let $F\subseteq O_K$ be
self-similar with respect to $\phi_i$ for $i=1$ to $n$. One can
define the fractal dimension of $F$ to be the real number $s$ for
which
$$
\sum_{i=1}^n Norm(a_i)^{-{s \over n_i}}=1.
$$
\begin{conj}
The above notion of fractal dimension for self-similar subsets of
$O_K$ is well-defined and well-behaved with respect to inclusion,
i.e. fractal dimension is independent of the choice of
self-similarities and compatible with inclusion of self-similar
subsets.
\end{conj}

\begin{rem}
 To stay on the safe side, it would have been better to state the general conjecture 
for linear polynomials were $n_i=1$ for all $i$, but we take this chance 
to state the more general conjecture, so that we indicate what we expect for the box formula of the fractal dimension.
\end{rem}

\subsection{Prospects of self-similarity in arithmetic geometry}

Extension of Siegel's theorem is proved using the following
strong version of Roth's theorem [Ra3]:

\begin{thm} (Fractal version of Roth's thereom on diopphantine approximation)
Fix a finitely generated field of characteristic zero $K$ 
and $\sigma :K\hookrightarrow \Cplx$ a
complex embedding. Let $V$ be a smooth projective algebraic
variety defined over $K$ and let $L$ be an very ample line-bundle on
$V$. Denote the arithmetic height function associated to the
line-bundle $L$ by $h_L$. Suppose $F\subset V(K)$ is a fractal
subset with respect to finitely many height-increasing
self-endomorphisms $\phi_i:V\to V$ defined over
$K$ such that for
all $i$ we have
$$
h_L(\phi_i(f))>m_ih_L(f)+0(1)+
$$
where $m_i>1$. Fix a Riemannian metric on $V_{\sigma}(\Cplx)$ and
let $d_{\sigma}$ denote the induced metric on
$V_{\sigma}(\Cplx)$. 
Then, for every $\delta>0$ and every choice
of an algebraic point $\alpha\in V(\bar {K})$ which is not a
critical value of any of the $\phi_i$'s and all choices of a
constant $C$, there are only finitely many fractal points
$\omega\in F$ approximating $\alpha$ such that 
$$
d_{\sigma}(\alpha ,\omega)\leq Ce^{-\delta h_L(\omega)}.
$$
\end{thm}

In fact an extended version of Siegel's theorem could be proved
using the above version of Roth's theorem [Ra3]:

\begin{thm}(Fractal version of Siegel's theorem on integral points) Fix a
finitely generated field of characteristic zero $K$. Let $V$ be a smooth affine algebraic variety
defined over $K$ with smooth projectivization $\bar V$ and let $L$
be an very ample line-bundle on $\bar V$. Denote the arithmetic height
function associated to the line-bundle $L$ by $h_L$. Suppose
$F\subset V({K})$ is a fractal subset with respect to
finitely many height-increasing polynomial self-endomorphisms
$\phi_i:V\to V$ defined over $K$ such that for all $i$ we have
$$
h_L(\phi_i(f)) > m_ih_L(f)+0(1)
$$
where $m_i>1$. One could also replace this assumption with norm
analogue. For any affine hyperbolic algebraic curve $X$ embedded
in $V$ defined over $K$ we have $X(K)\cap F$ is a
finite set.
\end{thm}

\section{self-similarity in category theory}

In category theory, one tries to translate properties of
mathematical objects by their external behavior, i.e. in terms of
morphisms to other objects. In some sense, category theory is
sociology of mathematical concepts, where dealing with internal
structures of mathematical objects is regarded as psychology of
these objects. It is believed that there is a duality between
sociological and psychological interpretations of mathematical
phenomena.

\subsection{Self-similar objects and self-similar morphisms}

Let us keep the example of self-similar vector spaces and
functions on Cantor fractal in mind. For objects $X$ and $Y$ in an
abelian category $\mathcal{C}$, one can restrict morphisms $X\to
Y$ to images of self-similarity maps $\phi_i:X\to X$. Therefore,
for Self-similar objects $X$, we get maps
$$
Hom(X,Y)\rightleftarrows Hom(\phi_i(X),Y)\cong Hom(X,Y)
$$
such that the natural induced maps
$$
\prod_i Hom(\phi_i(X),Y)\longrightarrow Hom(X,Y)\longrightarrow
\prod_i Hom(\phi_i(X),Y)
$$
are one-to-one correspondences which are inverse to each other.
There are also maps in the other direction:
$$
Hom(Y,X)\rightleftarrows Hom(Y,\phi_i(X))\cong Hom(Y,X)
$$
such that the natural maps
$$
\prod_i Hom(Y,\phi_i(X))\longrightarrow Hom(Y,X)\longrightarrow
\prod_i Hom(Y,\phi_i(X))
$$
are also one-to-one correspondences which are inverse to each
other. Although, it seems that working with abelian categories is
essential, one can axiomatize self-similarity in arbitrary
category as follows:

\begin{defn}
An object $X$ of a category $\mathcal{C}$ is called self-similar
if for any object $Y$ there are finitely many similarity maps
$$
\psi_i:Hom(X,Y)\to Hom(X,Y)
$$
such that the induced map
$$
\prod_i \psi_i:Hom(X,Y)\longrightarrow \prod_i Hom(X,Y)
$$
is a one-to-one correspondence. Similarly, there should exist
co-similarity maps
$$
\mu_j:Hom(Y,X)\to Hom(Y,X)
$$
such that the induced map
$$
\prod_j \mu_j:Hom(Y,X)\longrightarrow \prod_j Hom(Y,X)
$$
is a one-to-one correspondence.
\end{defn}

Morphisms between self-similar vector spaces motivate how to
define morphisms between self-similar objects in a category. but
we need an abelian category to have a one-to one correspondence
between similarity maps and co-similarity maps. In this case, for
self-similar objects $X$ and $Y$ with similarity maps indexed by
index sets $I$ and $J$ respectively, we have one-to-one
correspondences
$$
\prod_i Hom(X,Y)\longrightarrow Hom(X,Y)\longrightarrow \prod_j
Hom(X,Y).
$$
Suppose one fixes a surjective map $\pi:I\to J$. Elements of
$Hom(X,Y)$ are self-similar with respect to $\pi$ if the
combination of the above one-to-one correspondences respects the
indices according to $\pi$.
\subsection{Self-similar categories}

In order to formulate self-similar categories, the natural choice
for similarity maps would be endo-functors of the category. Here
is the first candidate coming to mind.

\begin{defn}
A category $\mathcal{C}$ is called a self-similar category with
respect to finitely many similarity functors
$\phi_i:\mathcal{C}\to \mathcal{C}$ if $\mathcal{C}$ is disjoint
union of its images under similarity maps:
$$
\mathcal{C}=\coprod_i \phi_i(\mathcal{C}).
$$
a self-similar morphism between self-similar categories is
defined similar to morphisms between self-similar vector spaces.
\end{defn}

On the other hand, a natural example for a self-similar category
should be the category of self-similar vector-spaces over a field
or the category of modules over a self-similar ring. In the first
case, there is no ambient concept of similarity maps. In the
second case, all similarity maps are surjective and thus they
have intersecting images.

\section{Some philosophical remarks}

Self-similar objects represent a form of co-finiteness which
embeds in the realm of finite mathematics because of their
recursive self-similarity. Box formulas let us think of
self-similar objects as finite dimensional structures, which
confirms finiteness of these structures. Parallel geometric,
algebraic and arithmetic formulations suggests that these objects
are fundamental in understanding mathematical structures and
glossary of examples show that they naturally appear in everyday
mathematics.

Fractal dimension of arithmetic self-similar objects encodes the
asymptotic behavior of points of bounded height. It does not
escape the eyes that for self-similar arithmetic objects, points
of bounded height are of polynomial growth. It would be
interesting if one could introduce an interpretation of the
coefficient of the asymptotic in terms of self-similarity
concepts.

Considering similarity as an abstract concept is confirmed to be
a natural abstraction. Because, it appears naturally in a number
of different mathematical frameworks and helps in generalizing and
unifying many seemingly unrelated results and conjectures. This
would be the main philosophical contribution of this manuscript.

Sharif University of Technology,  rastegar@sharif.ir

\end{document}